\documentclass[12pt]{article}
\usepackage[latin1]{inputenc}
\usepackage[T1]{fontenc}
\usepackage{amsmath,amssymb,amstext}
\usepackage{theorem}
\usepackage{a4wide}
\usepackage{multicol}
\usepackage[english]{babel}
\usepackage{enumerate}
\newcommand{\R}{\mathbb{R}}

\newcounter{moncompteur}
{}
\newtheorem{prop}[moncompteur]{\textbf{Proposition} }{}
\newtheorem{thm}[moncompteur]{\textbf{Theorem} }{}
\newtheorem{lemma}[moncompteur]{\textbf{Lemma} }{}
\newtheorem{defi}[moncompteur]{\textbf{Definition} }{}
\newtheorem{rk}[moncompteur]{\textbf{Remark} }{}
\newenvironment{dem}{\textbf{Proof}~:\\ }{\flushright$\blacksquare$\\}
\title{Harnesses, Lévy bridges and \textit{Monsieur Jourdain}}
\author{\textsc{Mansuy} Roger\footnote{contact~: \textbf{mansuy@ccr.jussieu.fr}}, \textsc{Yor} Marc}
\date{Laboratoire de Probabilit\'es et mod\`eles al\'eatoires,\\
Universit\'e Pierre
et Marie Curie-Paris VI,\\ Casier 188,\\ 4,place Jussieu,\\
F-75252 Paris Cedex 05}
\begin{document}
\maketitle \pagestyle{myheadings}
\begin{abstract}Relations between harnesses (\cite{MR36:7190},
\cite{williams}) and initial enlargements of the filtration of a
Lévy process with its positions at fixed times are investigated.
\end{abstract}
\noindent\textbf{Keywords~:} \emph{Harnesses, Lévy processes, past-future martingales,
enlargement of filtration}\\
\textbf{AMS 2000 subject classification ~:}*60G48, 60G51, 60G44,
60G10 \tableofcontents
\begin{section}{Introduction}
In order to model long-range misorientation within crystalline
structure of metals, Hammersley \cite{MR36:7190} introduced
various notions of processes which enjoy particular conditional
expectation properties. Among these, harnesses will be of
particular interest. Let us precise the definition~:
\begin{defi}:\\ Let $(H_t;t\geq 0)$ be a measurable process such that for all $t$, $\mathbb{E}[|H_t|]<\infty$, and
define for all $t<T$~:$$\mathcal{H}_{t,T}:=\sigma\left\{H_s; s\leq
t; H_u; u\geq T\right\}$$ $H$ is said to be a harness if, for all
$a<b<c<d$
\begin{eqnarray}\mathbb{E}\left[\frac{H_c-H_b}{c-b}|\mathcal{H}_{a,d}\right]&=&\frac{H_d-H_a}{d-a}
\label{eq5}\end{eqnarray}
\end{defi}
One might also define the notion of
$(\mathcal{F}_{t,T})_{t<T}$-harness as soon as
$\mathcal{H}_{t,T}\subset\mathcal{F}_{t,T}$, with obvious
hypothesis on a "past-future" filtration $\mathcal{F}$, which may
be just as useful as the notion of Brownian motion with respect to
a filtration. The equality may be reformulated as follows~: $H$ is
a harness if and only if for all $s<t<u$
\begin{eqnarray}\mathbb{E}\left[H_t|\mathcal{H}_{s,u}\right]&=&
\frac{t-s}{u-s}H_u+\frac{u-t}{u-s}H_s\label{def1}\end{eqnarray}
Such a formulation justifies that Harnesses are sometimes called
affine processes (See \cite{cy} chapter 6).\\
We note that Williams (\cite{MR50:15005} and \cite{williams})
proved the following striking result~:the only squared integrable
continuous harnesses are Brownian motions with drifts.
This latter result shows how rigid the property of being a
continuous harness is and may help understand why studies of
harnesses with continuous time were so few during the past twenty
years. On the other hand, some multi-parameter versions appeared,
imitating Williams arguments
(See \cite{MR82h:60152}, \cite{MR93h:60129}, \cite{MR91e:60158} and \cite{MR999888}).\\
Glancing through the literature, it seems that no study of
discontinuous harnesses has been performed. Our reference to
Monsieur Jourdain (a character of Molière (1622-1673)
\cite{moliere}) in the title  alludes to this point; as Monsieur
Jourdain discovers he was practising prose without being aware of
it, the following theorem shows that a number of authors have been
dealing with harnesses.
\begin{thm}:
\begin{enumerate}[(i)]
  \item (Jacod-Protter, \cite{MR89j:60102})Let $(\xi_t; t\geq
  0)$ be an integrable Lévy process (that is: $\forall t,
  \mathbb{E}[|\xi_t|]<\infty$) and define
  $$\mathcal{F}_{t,T}=\sigma\left\{\xi_s;s\leq t;\xi_u;u\geq
  T\right\}$$ Then for any given $T>0$, there is the
  decomposition formula~:\begin{eqnarray}
\xi_t&=&M_t^{(T)}+\int_0^tds\ \frac{\xi_T-\xi_s}{T-s}
\label{eq4}\end{eqnarray}where $(M_t;t\leq T)$ is a
$\left(\mathcal{F}_{t,T};t\leq T\right)$-martingale
  \item In a general framework, an integrable process $(H_t;t\geq
  0)$ is a $\left(\mathcal{F}_{t,T}\right)$-harness if and only
  if, for every $T>0$, there exists $(M_t^{(T)})_{t<T}$ a $\left(\mathcal{F}_{t,T}
  ;t< T\right)$-martingale such that\begin{eqnarray}\forall t< T,\
H_t&=&M_t^{(T)}+\int_0^tds\ \frac{H_T-H_s}{T-s}\label{eq6}
\end{eqnarray}
\end{enumerate}
\end{thm}
For further results along this line, see Exercise 6.19 in
\cite{cy} which provides a few references about harnesses. In the
particular case of a Brownian motion $\xi$, formula (\ref{eq4})
may be attributed to Itô \cite{MR81f:60082} but was already
sketched by Lévy \cite{MR6:5e} and \cite{MR7:314g}. See also
Jeulin-Yor \cite{MR82d:60084}. Our motivation for writing this
note is that harnesses -through formula (\ref{eq4})- seem to
become more topical; indeed some recent works (\cite{oksendal} and
\cite{kh}) develop financial models of markets with well informed
agents (also called insiders) where formula (\ref{eq6}) plays a
key-role. Some other papers (\cite{math.PR/0210218} or
\cite{math.PR/0312402}) also deal with some notions of harness
derived directly from the pioneering work of Hammersley,
but are apparently far from the preceding discussion.\\
This note is organized as follows~:
\begin{itemize}
  \item First we prove part (ii) of the theorem.
  \item Section \ref{gir} is devoted to an alternative proof of
  the decomposition formula (\ref{eq4}) of Jacod-Protter \cite{MR89j:60102}
  thanks to the absolute continuity of the law of a Lévy process
  and its bridge.
  \item In Section \ref{paf}, we develop the more general notion
  of past-future martingale and provide as many examples as
  possible.
\end{itemize}
\end{section}
\begin{section}{Relations between Lévy bridges and
harnesses}\label{lh}
\begin{paragraph}{(\ref{lh}.1)}Let $(B_t;t\geq 0)$ be a 1-dimensional Brownian
motion; it is well known that a realization of the Brownian bridge
over the time interval $[0,T]$, starting at $x$ and ending at $y$,
is:
\begin{eqnarray}&\left\{x+\left(B_t-\frac{t}{T}B_T\right)+\frac{t}{T}y;
t\leq T\right\}&\label{eq1}\end{eqnarray} Moreover, the
semimartingale decomposition of this bridge is also well-known; it
is the solution of the SDE~:
\begin{eqnarray}X_t&=&x+\beta_t+\int_0^tds\ \frac{y-X_s}{T-s};\ \
t\leq T \label{eq2}
\end{eqnarray}
where $(\beta_t;t\leq T)$ is a standard Brownian motion.\\
This decomposition formula (\ref{eq2}) is, in fact, equivalent to
the semimartingale decomposition of $(B_t;t\leq T)$ in the
enlarged filtration
$\mathcal{B}_t^{(T)}:=\mathcal{B}_t\vee\sigma(B_T)$, where
$\mathcal{B}_t=\sigma\{B_s; s\leq t\}$~:
\begin{eqnarray}
B_t&=&\gamma_t^{(T)}+\int_0^tds\ \frac{B_T-B_s}{T-s}\label{eq3}
\end{eqnarray}
where $(\gamma_t^{(T)};t\leq T)$is a $(\mathcal{B}_t^{(T)}; t\leq
T)$-Brownian motion; in particular, it is independent of $B_T$.
See \cite{MR81f:60082} and \cite{MR82d:60084} for a discussion of
(\ref{eq2}) and (\ref{eq3}).
\end{paragraph}
\begin{paragraph}{(\ref{lh}.2)} It has been shown by Jacod-Protter
\cite{MR89j:60102} that formula (\ref{eq3}) in fact extends to any
integrable Lévy process $(\xi_t; t\geq 0)$ in the following way~:
\begin{eqnarray}
\xi_t&=&M_t^{(T)}+\int_0^tds\ \frac{\xi_T-\xi_s}{T-s}
\end{eqnarray}where $(M_t^{(T)};t\leq T)$ is a martingale in
the enlarged filtration
$\mathcal{F}^{(T)}_t=\mathcal{F}_t\vee\sigma(\xi_T)$, where
$\mathcal{F}_t=\sigma(\xi_s; s\leq t)$.
\end{paragraph}
\begin{paragraph}{(2.3)} Here is the proof of part (ii) of Theorem 2:
\begin{enumerate}[a.]
    \item $(\Rightarrow)$ Let $H$ be a harness and $s<t<T$.\\
    Define $M_t^{(T)}=H_t-\int\limits_0^t\frac{H_T-H_u}{T-u}du$.\\
    Then, the harness property implies
    \begin{eqnarray*}
    \mathbb{E}\left[M_t^{(T)}|\mathcal{F}_{s,T}\right]&=&
    \mathbb{E}\left[H_t|\mathcal{F}_{s,T}\right]-\int\limits_0^s\frac{H_T-H_u}{T-u}du-\int\limits_s^t\mathbb{E}\left[\frac{H_T-H_u}{T-u}
    |\mathcal{F}_{s,T}\right]du\\
    &=&\frac{T-t}{T-s}H_s+\frac{t-s}{T-s}H_T-\int\limits_0^s\frac{H_T-H_u}{T-u}du-\int\limits_s^t\frac{H_T-H_s}{T-s}du\\
    &=&M_s^{(T)}
    \end{eqnarray*}
    \item $(\Leftarrow)$ First, remark it is enough to show that,
    for all $s<t<T$
    \begin{eqnarray}
    \mathbb{E}\left[\frac{H_t-H_s}{t-s}|\mathcal{F}_{s,T}\right]&=&\frac{H_T-H_s}{T-s}\label{eq7}\end{eqnarray}
    Indeed, if $r<s<t<T$, then
    \begin{eqnarray*}
    \mathbb{E}\left[\frac{H_t-H_s}{t-s}|\mathcal{F}_{r,T}\right]&=&\mathbb{E}\left[\mathbb{E}\left[\frac{H_t-H_s}{t-s}
    |\mathcal{F}_{s,T}\right]|\mathcal{F}_{r,T}\right]\\
    &=&\mathbb{E}\left[\frac{H_T-H_s}{T-s}|\mathcal{F}_{r,T}\right]\\
    &=&\frac{\mathbb{E}\left[H_r-H_s|\mathcal{F}_{r,T}\right]}{T-s}+\frac{H_T-H_r}{T-s}\\
    &=&\frac{r-s}{T-s}\frac{H_T-H_r}{T-r}+\frac{H_T-H_r}{T-s}\\
    &=&\frac{H_T-H_r}{T-r}
    \end{eqnarray*}
    It only remains to prove formula (\ref{eq7}). The
    assumed decomposition formula (\ref{eq6}) yields to
    \begin{eqnarray*}
    H_t-H_s&=&M_t^{(T)}-M_s^{(T)}+\int\limits_s^tdv\frac{H_T-H_v}{T-v}\\
    \end{eqnarray*}
    Therefore
    \begin{eqnarray*}
    \mathbb{E}\left[H_t-H_s|\mathcal{F}_{s,T}\right]&=&
    \int\limits_s^tdv\frac{\mathbb{E}\left[H_T-H_v|\mathcal{F}_{s,T}\right]}{T-v}\\
    &=& \int\limits_s^t\frac{dv}{T-v}\left(H_T-H_s\right)-\int\limits_s^t\frac{dv}{T-v}
    \mathbb{E}\left[H_v-H_s|\mathcal{F}_{s,T}\right]\\
    \end{eqnarray*}
    Hence, $s$ and $T$ being fixed,
    $\phi(t):=\mathbb{E}\left[H_t-H_s|\mathcal{H}_{s,T}\right]$ solves the following first
    order linear differential equation~:
    $$\phi(t)=\int_s^t\frac{dv}{T-v}\left(H_T-H_s\right)-\int\limits_s^t\frac{dv}{T-v}\phi(v);\ s\leq t\leq T$$
    But this equation admits only one solution vanishing at $s$ and a standard computation yields to
    $\phi(t)=\frac{H_T-H_s}{T-s}(t-s)$ which is formula
    (\ref{eq7}).
\end{enumerate}
\begin{rk}:\\Contrary to the very definition of harness, this
proposition exhibits a privileged direction of time. So a similar
representation property with the opposite time-direction can be
derived. Namely, a measurable process $H$ is a harness on $[0,T]$,
if and only if, for all $T>\tau>0$, there exists
$(N_t^{(\tau)};\tau< t\leq T)$ a $(\mathcal{F}_{\tau,t}; \tau<
t\leq T)$-reverse martingale such that
\begin{eqnarray}\forall \tau<t\leq T,\
H_t&=&N_t^{(\tau)}-\int_t^Tds\frac{H_{\tau}-H_s}{\tau-s}
\end{eqnarray}
\end{rk}
\end{paragraph}
\end{section}
\begin{section}{A Girsanov proof of the decomposition
formula}\label{gir}
\begin{paragraph}{(\ref{gir}.1)}It is well known (see e.g. \cite{MR95i:60070}) that the law of the
bridge of a Markov process is locally equivalent to the law of the
"good" Markov process, more precisely, if $X$ is a Markov process
with $p_t(x,y)$ as its semigroup density from $x$ to $y$, then the
following absolute continuity relationship between
$\mathbb{P}^t_{x\rightarrow y}$, the law of the bridge of length
$t$ from $x$ to $y$ and $\mathbb{P}_x$ the law of $X$ starting at
$x$ holds~: \begin{eqnarray}\mathbb{P}_{x\rightarrow
y|\mathcal{F}_s}^t&=&\frac{p_{t-s}(X_s,y)}{p_t(x,y)}.\mathbb{P}_{x|\mathcal{F}_s}\label{mark}\end{eqnarray}
If $\xi$ is a Lévy process, $\phi_t(\centerdot)$ will denote the
density of the law of $\xi_t$, assuming it exists (see
\cite{MR2003b:60064}
 for conditions on a Lévy process to have such a density). The equality (\ref{mark}) then becomes~:
\begin{eqnarray}\mathbb{P}_{x\rightarrow
y|\mathcal{F}_s}^t&=&\frac{\phi_{t-s}(y-\xi_s)}{\phi_t(y-x)}.\mathbb{P}_{x|\mathcal{F}_s}
\label{levyabs}\end{eqnarray} We now stay in the context of a Lévy
process.
\begin{lemma}:\\ If $(M^y_t;t\leq T,y\in\R)$ denote a family of variables such that
\begin{itemize} \item for any $y\in\R$, $(M^y_t;t\leq T)$ is  a $P^T_{x\rightarrow
y}$-martingale. \item $(t,y)\mapsto M_t^y$ is
measurable.\end{itemize}Then $(M^{\xi_T}_t;t\leq T)$ remains a
$P_x$-martingale with respect to the filtration initially enlarged
with $\xi_T$.
\end{lemma}
\begin{dem}Let $(M^y_t;t\leq T,y\in\R)$ be such a family of $P^T_{x\rightarrow y}$-martingales; then,  for all $s<t<T$
and $\Gamma_s\in\sigma(\xi_u;u\leq s )$,
\begin{eqnarray*}\mathbb{E}_{x\rightarrow
y}^{T}\left[1_{\Gamma_s}(M_t^{y}-M_s^{y})\right]&=&0\end{eqnarray*}
This implies, for any  bounded Borel function $f$,
\begin{eqnarray*}\int \mathbb{P}_x (\xi_T\in
dy)f(y)\mathbb{E}_{x\rightarrow
y}^{T}\left[1_{\Gamma_s}(M_t^{y}-M_s^{y})\right]&=&0\end{eqnarray*}
Therefore
\begin{eqnarray*}
\mathbb{E}_x\left[f(\xi_T)1_{\Gamma_s}(M_t^{\xi_T}-M_s^{\xi_T})\right]&=&0
\end{eqnarray*}
So, $M_t^{\xi_T}$ a $P_x$-martingale with respect to the
filtration enlarged with $\xi_T$.
\end{dem}
\end{paragraph}
\begin{paragraph}{(\ref{gir}.2)}If we suppose, without any loss of generality, that $\mathbb{E}[\xi_1]=0$,
then $\xi$ is a $P_x$-martingale (in any other case, we will study
the Lévy process $\xi_t-\mbox{\texttt{d}}t$ where \texttt{d} is
the drift term of $\xi$). We shall denote $(\sigma^2,\nu)$ its
local characteristics (Brownian term and Lévy measure) and
$\mathcal{L}$ its infinitesimal generator. For the sake of
simplicity, note that $\tilde{\mathcal{L}}$, the infinitesimal
generator of the time-space process $(t,\xi_t)$ satisfies
$$\tilde{\mathcal{L}}=\frac{\partial\ }{\partial t}+\mathcal{L}$$
Thanks to the Girsanov theorem and the absolute continuity
relationship (\ref{levyabs}), the process
$$\xi_t-\int_0^t\frac{d\langle
\xi_\centerdot,\phi_{T-\centerdot}(y-\xi_\centerdot)\rangle_s}{\phi_{T-s}(y-\xi_s)}$$
defines a $P^T_{x\rightarrow y}$-martingale and therefore
$$\xi_t-\int_0^t\frac{d\langle
\xi_\centerdot,\phi_{T-\centerdot}(\xi_T-\xi_\centerdot)\rangle_s}{\phi_{T-s}(\xi_T-\xi_s)}$$
is a $P_x$-martingale with respect to the filtration enlarged with
$\xi_T$; this process will now be compared with
$(M^{(T)}_t)_{t\leq T}$ in part (ii) of Theorem 2. Namely, we aim
to prove that
\begin{eqnarray}
\langle
\xi_\centerdot,\phi_{T-\centerdot}(y-\xi_\centerdot)\rangle_t&=&
\int\limits_0^t\frac{y-\xi_s}{T-s}\ \phi_{T-s}(y-\xi_s)ds
\end{eqnarray}
that is, with our notation~:
\begin{eqnarray*}
\tilde{\mathcal{L}}(x\phi_{T-s}(y-x))(s,x)&=&\frac{y-x}{T-s}\phi_{T-s}(y-x)
\end{eqnarray*}
Now,
\begin{eqnarray*}
\tilde{\mathcal{L}}(x\phi_{T-s}(y-x))(s,x)&=&-\sigma^2\phi'_{T-s}(y-x)
+\int\nu(dz)z\phi_{T-s}(y-x-z)
\end{eqnarray*}[This
computation is quite easy once we note that $(t,x)\mapsto
\phi_{T-t}(y-x)$ is a space-time harmonic function.] \\The
following lemma concludes the proof~:
\begin{lemma}:\\For any integrable Lévy process with local characteristics
$(\sigma^2,\nu)$ and transition probability density $\phi$,
\begin{eqnarray}
-\sigma^2\phi'_{u}(x)
+\int\nu(dz)z\phi_{u}(x-z)&=&\frac{x}{u}\phi_{u}(x)\label{machin}\end{eqnarray}
\end{lemma}
\begin{dem}
From the very definition of the Lévy exponent, we have~:
\begin{eqnarray}
\int e^{i\lambda x}\phi_u(x)dx=\mathbb{E}\left[e^{i\lambda \xi_u}
\right]&=& e^{-u\Phi(\lambda)}\label{simple}
\end{eqnarray}
Differentiation in $\lambda$ within this equality yields to
\begin{eqnarray*}
i\int  x\phi_u(x)e^{i\lambda x}dx&=&-u\Phi'(\lambda)
e^{-u\Phi(\lambda)}
\end{eqnarray*}
with $\Phi'(\lambda)=\sigma^2\lambda-i\int \nu(dz)ze^{i\lambda
z}$\\ Replacing $e^{-u\Phi(\lambda)}$ with the expression in
(\ref{simple}) and noting that
\begin{eqnarray*}
\lambda\int \phi_u(x)e^{i\lambda x}dx&=&i\int \phi'_u(x)e^{i\lambda x}dx\\
\int \nu(dz)ze^{i\lambda z}\int \phi_u(x)e^{i\lambda x}dx &=& \int
dx e^{i\lambda x}\int \nu(dz)z\phi_u(x-z)
\end{eqnarray*}
we obtain~:
\begin{eqnarray*}
i\int  x\phi_u(x)e^{i\lambda x}dx&=&-u\int dx \left(-\sigma^2
\phi_u'(x)+\int\nu(dz)z\phi_u(x-z) \right)e^{i\lambda x}
\end{eqnarray*}
\end{dem}
\begin{rk}:\\The right-hand side of (\ref{machin}) can also be interpreted, for
skip-free Lévy processes, as the density of the first hitting time
thanks  to Kendall's identity (See e.g.
\cite{MR2002i:60099}).\end{rk}
\end{paragraph}
\end{section}
\begin{section}{A wider class of processes: the past-future
martingales}\label{paf}
\begin{paragraph}{(\ref{paf}.1)} If $\mathcal{F}$ denotes a past-future filtration,
 the following definition generalizes the
notion of a $\mathcal{F}$-harness~:
\begin{defi}:\\The two-parameters process $(M_{s,t})_{0\leq
s<t<\infty}$ is said to be a past-future martingale with respect
to $(\mathcal{F}_{s,t})_{0\leq s<t<\infty}$ if~:
\begin{enumerate}
  \item $\forall s<t$, $\mathbb{E}\left[|M_{s,t}|\right]<\infty$
  \item $\forall s<t$, $M_{s,t}$ is
  $\mathcal{F}_{s,t}$-measurable.
  \item $\forall r<s<t<u$,
  $\mathbb{E}\left[M_{s,t}|\mathcal{F}_{r,u}\right]=M_{r,u}$
\end{enumerate}
\end{defi}
\begin{rk}:
\begin{itemize}
  \item As previously mentioned, a process $H$ is a $\mathcal{F}$-harness if and
  only if $\left(\frac{H_t-H_s}{t-s}\right)_{0\leq s<t<\infty}$ is
  a past-future martingale.
  \item Note that past-future
  martingales are reverse martingales indexed by the intervals of
  $\R^+$.
\end{itemize}
\end{rk}
\end{paragraph}
\begin{paragraph}{(\ref{paf}.2)} Here we are to detail some non trivial
past-future martingales related to a standard Brownian motion
$(B_t;t\geq 0)$.
\begin{enumerate}
    \item Let $f_+$ and $f_-$ be two both square-integrable and
    integrable
    functions on $\R^+$ and $C\in\R$. Then the process
    $(M_{s,t})_{0\leq s<t<\infty}$ defined for all $s<t$ by~:
    \begin{eqnarray*}
    M_{s,t}&=&\int_0^sf_-(u)dB_u+\int_t^{\infty}f_+(u)dB_u+\ldots\\
    &&\ldots+
    \frac{B_t-B_s}{t-s}\left(C-\int_0^sf_-(u)du-\int_t^{\infty}f_+(u)du\right)
    \end{eqnarray*} is a past-future Brownian martingale.\\ One notices that the stochastic integral terms
    associated to the functions $f_{\pm}$ have to be "compensated"
    with a harness term.
    \item An exponential example can easily be derived from this
    latter. Within the same framework, the two-parameter process $(N_{s,t})_{0\leq s<t\infty}$
     defined for all $s<t$\begin{eqnarray*}\ln
    N_{s,t}&=&M_{s,t}+\cfrac12\int\limits_0^sf_-^2(u)du
    +\cfrac12\int\limits_t^{\infty}f_+^2(u)du+\ldots\\&&\ldots+\cfrac{t-s}2
\left(C-\int\limits_0^sf_-(u)du-\int\limits_t^{\infty}f_+(u)du\right)^2\end{eqnarray*}
is a past-future martingale.
\end{enumerate}
\end{paragraph}
\begin{paragraph}{(\ref{paf}.3)}Previous examples can easily be
extended to more general Lévy processes~:
\begin{prop}:\\ Let $\xi$ be a Lévy process and $f$ an
integrable function with locally finite variation, chosen to be
right-continuous with left limits, such that $\int_0^{\infty}
f(u-)d\xi_u$ exists. Then, for all $s<t$,
\begin{eqnarray*}
M_{s,t}&=&\int_0^sf(u^-)d\xi_u+\int_t^{\infty}f(u^-)d\xi_u
+\frac{\xi_t-\xi_s}{t-s}\int_s^tf(u)du+[\xi_\centerdot,f(\centerdot)]_s
-[\xi_\centerdot,f(\centerdot)]_t
\end{eqnarray*}
defines a past-future martingale.
\end{prop}
\begin{dem}Indeed thanks to integration by parts formula
\begin{eqnarray*}
\mathbb{E}\left[\int_s^tf(u^-)d\xi_u|\xi_t,\xi_s\right]&=&f(t^-)\xi_t-f(s)\xi_s-
\int_s^t\mathbb{E}\left[\xi_{u^-}|\xi_t,\xi_s\right]df(u)+[\xi_\centerdot,f(\centerdot)]_s
-[\xi_\centerdot,f(\centerdot)]_t\\
&=&\xi_t\left(f(t^-)-\int_s^t\frac{u-s}{t-s}df(u)\right)
+\xi_s\left(-f(s)-\int_s^t\frac{t-u}{t-s}df(u)\right)\\
&&\ldots+[\xi_{\centerdot},f(\centerdot)]_s
-[\xi_{\centerdot},f(\centerdot)]_t\\
&=&\frac{(\xi_t-\xi_s)}{t-s}\int_s^tf(u)du+[\xi_\centerdot,f(\centerdot)]_s
-[\xi_\centerdot,f(\centerdot)]_t
\end{eqnarray*}
Therefore
\begin{eqnarray*}
M_{s,t}&=&\mathbb{E}\left[\int_0^{\infty}f(u^-)d\xi_u|\mathcal{H}_{s,t}\right]
\end{eqnarray*}
\end{dem}
\end{paragraph}
\end{section}
\nocite{*}
\bibliographystyle{alpha}
\bibliography{harnais}

\end{document}